\documentclass[12pt,twoside]{article}

\usepackage{geometry} 
\usepackage{graphicx} 

\usepackage{booktabs} 
\usepackage{array} 
\usepackage{paralist} 
\usepackage{verbatim} 
\usepackage{subfig} 

\usepackage{amssymb}
\usepackage{amsmath}
\usepackage{fancyref}
\usepackage{amsthm}
\usepackage{cases}
\usepackage{hyperref} 

\usepackage{fancyhdr} 
\usepackage{sectsty}
\allsectionsfont{\sffamily\mdseries\upshape} 

\usepackage[nottoc,notlof,notlot]{tocbibind} 
\usepackage[titles,subfigure]{tocloft} 

\numberwithin{equation}{section}

\title{\Large{\bf General Methods For Solving Ordinary Differential Equations 1}}
\author{\bf {Daiyuan Zhang}}
\begin{document}
\maketitle
\centerline{College of Computer}
\centerline{Nanjing University of Posts and Telecommunications}
\centerline{Nanjing, P.R. China}
\centerline{dyzhang@njupt.edu.cn, zhangdaiyuan2011@sina.com}

\newtheorem{Theorem}{\quad Theorem}[section]
\newtheorem{Definition}[Theorem]{\quad Definition}
\newtheorem{Corollary}[Theorem]{\quad Corollary}
\newtheorem{Lemma}[Theorem]{\quad Lemma}
\newtheorem{Example}[Theorem]{\quad Example}

\centerline{}

\begin{abstract}
	All the theories and methods proposed in this paper are my original innovations. On the basis of the new mathematical expansion I proposed in 2014, this paper presents a general method for solving linear ordinary differential equations. Different from traditional known methods, this paper deals with a completely new method for solving linear ordinary differential equations. It is concluded that, under some conditions which easily be satisfied in practical applications, any order of linear ordinary differential equations can be solved by the method of recursion and reduction of order step by step. 
\end{abstract}


{\footnotesize Copyright $\copyright$ This is my original work. Prohibit any plagiarism.}

\section{Introduction}\label{sec1}

Based on the new mathematical expansion I proposed in \cite{daiyuan_zhang_new_2014} (Daiyuan Zhang, 2014),  this paper studies a new method for solving the $k$th order linear ordinary differential equations, where $k$ is a positive integer. 

The mainly known methods for solving the $k$th order linear ordinary differential equations can be briefly discussed as follows. 

The solution to the first order linear ordinary differential equation has already been mature, and people can get its general solution very easily.

The solution of the second order linear ordinary differential equation is not so easy. Although one solution can be obtained when another solution is known, if all of the solutions are unknown, it will bring difficulties.

It is even more difficult to solve linear ordinary differential equations when $k\ge3$.

For a given linear ordinary differential equations with constant coefficients, the method of characteristic equation can be used to solve it, but it is not suitable for solving linear ordinary differential equations with variable coefficients.

Generally speaking, although the infinite series method can get some solutions, it is not easy to determine the convergence of the infinite series. Therefore, it is difficult to use the initial conditions to determine the related coefficients.

Although the numerical methods can obtain some practical solutions, and can also obtain numerical solutions with high accuracy in some cases, the numerical method is difficult to see the inherent law for solving the differential equations.

For the sake of convenience, we call those known methods mentioned above as traditional methods in this paper. 

The traditional methods for solving differential equations can be found in many documents. For example, beginners of studying differential equations can read some introductory-level references such as \cite{william_a._adkins_ordinary_2010}, \cite{shair_ahmad_textbook_2014}, \cite{mircea_v._soare_ordinary_2007} etc. More advanced topics about differential equations, please refer to literature \cite{jack_k._hale_ordinary_1980} etc. 

In this paper, I present a unified method for solving linear ordinary differential equations. The original foundation of this paper is the new mathematical expansion proposed in my paper \cite{daiyuan_zhang_new_2014}, without this new mathematical expansion \cite{daiyuan_zhang_new_2014}, we can not get the new method for solving the high order linear ordinary differential equations. This method is completely different from the known traditional methods. This method can not only solve linear ordinary differential equations with constant coefficients, but also solve linear ordinary differential equations with variable coefficients. It can not only solve low order linear ordinary differential equations, but also solve high order linear ordinary differential equations. The procedure for solving linear ordinary differential equations proposed in this paper is reduction of order and recursive calculation. For a given linear ordinary differential equation, the order of the linear ordinary differential equation can be gradually reduced by adding, subtraction, multiplication, division and derivative operation from the coefficients of the given equation, and then the solution can be obtained at the end, and the  conditions are automatically included in the solution.

\section{Theoretical basis }\label{sec2}

The research results in this section establish solid theoretical and practical basis for the new achievements proposed in this paper. This section contains two subsections. Subsection \ref{subsec2.1} gives a new proof of the new mathematical expansion proposed in paper \cite{daiyuan_zhang_new_2014}. Compared with document \cite{daiyuan_zhang_new_2014}, the new proof is simpler. In subsection \ref{subsec2.2}, based on my new mathematical expansion \cite{daiyuan_zhang_new_2014}, two important theorems are proved, they are the foundations for constructing differential equations, which can be used for solving the $k$th order linear ordinary differential equations (see section\ref{sec3}).

\subsection{A new mathematical expansion}\label{subsec2.1}

An important theorem concerning a new mathematical expansion for a given function $f\left( x \right)$ has been proposed in my paper \cite{daiyuan_zhang_new_2014} and given in the following:
\begin{Theorem}[See \cite{daiyuan_zhang_new_2014}]\label{thm2.1} 
Let $\mathbf{I}=\left[ A,\ B \right]$, $f\left( x \right)$  be a function having  finite $ (n+1) $th derivative ${f^{\left( n+1 \right)}}\left( x \right)$  everywhere in an open interval $\left( {A,\;B} \right)$  and assume that $f\left( x \right) \in {C^n}\left[ {A,\;B} \right]$, $a \in \left( {A,\;B} \right)$, then, for every $ x $ in $\left( {A,\;B} \right)$, we have
\begin{equation}\label{eq2.1}
	\begin{aligned}
		f\left( x \right) = {Z_n}\left[ {f\left( {x;\;a} \right)} \right] + {R_n}\left[ {f\left( {x;\;a} \right)} \right]
	\end{aligned}
\end{equation}
where
\begin{equation}\label{eq2.2}
	\begin{aligned}
		{Z_n}\left[ {f\left( {x;\;a} \right)} \right] = f\left( a \right) + \sum\limits_{k = 1}^n {{{\left( { - 1} \right)}^{k - 1}}\frac{{{f^{\left( k \right)}}\left( x \right)}}{{k!}}{{\left( {x - a} \right)}^k}}
	\end{aligned}
\end{equation}
\begin{equation}\label{eq2.3}
	\begin{aligned}
		{R_n}\left[ {f\left( {x;\;a} \right)} \right] = \frac{1}{{n!}}\int\limits_a^x {{{\left( {a - t} \right)}^n}{f^{\left( {n + 1} \right)}}\left( t \right)dt} 
	\end{aligned}
\end{equation}
where ${Z_n}\left[ {f\left( {x;\;a} \right)} \right]$ is known as the expansion of function $f\left( x \right)$  up to the $ n $-th derivative at the point $x = a$, and ${R_n}\left[ {f\left( {x;\;a} \right)} \right]$  is known as the remainder of the expansion ${Z_n}\left[ {f\left( {x;\;a} \right)} \right]$  at the point $x = a$.
\end{Theorem}
Where ${{Z}_{n}}\left[ f\left( x;\ a \right) \right]$ is also denoted as ${{Z}_{n}}\left[ f \right]$, and ${{R}_{n}}\left[ f\left( x;\ a \right) \right]$ is also denoted as ${{R}_{n}}\left[ f \right]$.

The proof of the above theorem was given in my paper \cite{daiyuan_zhang_new_2014}. For the completeness of this paper, a simplified proof of theorem \ref{thm2.1} is discussed below. 
\begin{proof}[Simplified Proof of Theorem \ref{thm2.1}]
Using the variable upper limit integral formula in Calculus, we have
\begin{equation}\label{eq2.4}
	\frac{d}{dx}\left( \frac{1}{n!}\int\limits_{a}^{x}{{{\left( a-t \right)}^{n}}{{f}^{\left( n+1 \right)}}\left( t \right)dt} \right)\text{=}{{\left( -1 \right)}^{n}}\frac{{{f}^{\left( n+1 \right)}}\left( x \right)}{n!}{{\left( x-a \right)}^{n}}
\end{equation}
i.e.
\begin{equation}\label{eq2.5}
	\begin{aligned}
		{f}'\left( x \right)& ={f}'\left( x \right)-{{\left( -1 \right)}^{n}}\frac{{{f}^{\left( n+1 \right)}}\left( x \right)}{n!}{{\left( x-a \right)}^{n}} +\frac{d}{dx}\left( \frac{1}{n!}\int\limits_{a}^{x}{{{\left( a-t \right)}^{n}}{{f}^{\left( n+1 \right)}}\left( t \right)dt} \right)  
	\end{aligned}
\end{equation}
whlie
\begin{equation}\label{eq2.6}
	\begin{aligned}
		& {f}'\left( x \right)-{{\left( -1 \right)}^{n}}\frac{{{f}^{\left( n+1 \right)}}\left( x \right)}{n!}{{\left( x-a \right)}^{n}} \\ 
		& ={f}'\left( x \right)+\sum\limits_{k=1}^{n-1}{{{\left( -1 \right)}^{k}}\frac{{{f}^{\left( k+1 \right)}}\left( x \right)}{k!}{{\left( x-a \right)}^{k}}}-\sum\limits_{k=1}^{n}{{{\left( -1 \right)}^{k}}\frac{{{f}^{\left( k+1 \right)}}\left( x \right)}{k!}{{\left( x-a \right)}^{k}}} \\ 
		& ={f}'\left( x \right)+\sum\limits_{k=2}^{n}{{{\left( -1 \right)}^{k-1}}\frac{{{f}^{\left( k \right)}}\left( x \right)}{\left( k-1 \right)!}{{\left( x-a \right)}^{k-1}}}+\sum\limits_{k=1}^{n}{{{\left( -1 \right)}^{k-1}}\frac{{{f}^{\left( k+1 \right)}}\left( x \right)}{k!}{{\left( x-a \right)}^{k}}} \\ 
		& =\sum\limits_{k=1}^{n}{{{\left( -1 \right)}^{k-1}}\left( \frac{{{f}^{\left( k \right)}}\left( x \right)}{\left( k-1 \right)!}{{\left( x-a \right)}^{k-1}} \right)}+\sum\limits_{k=1}^{n}{{{\left( -1 \right)}^{k-1}}\frac{{{f}^{\left( k+1 \right)}}\left( x \right)}{k!}{{\left( x-a \right)}^{k}}} \\ 
		& =\sum\limits_{k=1}^{n}{{{\left( -1 \right)}^{k-1}}\left( \frac{{{f}^{\left( k \right)}}\left( x \right)}{\left( k-1 \right)!}{{\left( x-a \right)}^{k-1}}+\frac{{{f}^{\left( k+1 \right)}}\left( x \right)}{k!}{{\left( x-a \right)}^{k}} \right)} \\ 
		& \text{=}\frac{d}{dx}\left( \sum\limits_{k=1}^{n}{{{\left( -1 \right)}^{k-1}}\frac{{{f}^{\left( k \right)}}\left( x \right)}{k!}{{\left( x-a \right)}^{k}}} \right) \\ 
	\end{aligned}	
\end{equation}

Substituting equation \eqref{eq2.6} into \eqref{eq2.5}, we have
\begin{equation}\label{eq2.7}
	{f}'\left( x \right)=\frac{d}{dx}\left( \sum\limits_{k=1}^{n}{{{\left( -1 \right)}^{k-1}}\frac{{{f}^{\left( k \right)}}\left( x \right)}{k!}{{\left( x-a \right)}^{k}}}+\frac{1}{n!}\int\limits_{a}^{x}{{{\left( a-t \right)}^{n}}{{f}^{\left( n+1 \right)}}\left( t \right)dt} \right)
\end{equation}

Therefore
\begin{equation}\label{eq2.8}
	f\left( x \right)=\sum\limits_{k=1}^{n}{{{\left( -1 \right)}^{k-1}}\frac{{{f}^{\left( k \right)}}\left( x \right)}{k!}{{\left( x-a \right)}^{k}}}+\frac{1}{n!}\int\limits_{a}^{x}{{{\left( a-t \right)}^{n}}{{f}^{\left( n+1 \right)}}\left( t \right)dt}+C
\end{equation}

In the above formula, $C$ is an arbitrary constant, or $C$ is an arbitrary amount independent of $x$. If $x=a$, $f\left( x \right)=f\left( a \right)$, then from formula \eqref{eq2.8}, we have
\begin{equation}\label{eq2.9}
	f\left( x \right)=f\left( a \right)+\sum\limits_{k=1}^{n}{{{\left( -1 \right)}^{k-1}}\frac{{{f}^{\left( k \right)}}\left( x \right)}{k!}{{\left( x-a \right)}^{k}}}+\frac{1}{n!}\int\limits_{a}^{x}{{{\left( a-t \right)}^{n}}{{f}^{\left( n+1 \right)}}\left( t \right)dt}
\end{equation}
\end{proof}

Formula \eqref{eq2.9} is the same as that in literature \cite{daiyuan_zhang_new_2014}. 

It must be pointed out that formula \eqref{eq2.9} is quite different from the Taylor expansion of the function $f\left( x \right)$. Formula \eqref{eq2.9} shows that the function $f\left( x \right)$ can be expressed by the derivatives ${{f}^{\left( k \right)}}\left( x \right)$ of the function $f\left( x \right)$, but the Taylor expansion of function $f\left( x \right)$ is a polynomial.

\subsection{Using the new mathematical expansion to construct differential equations}\label{subsec2.2}

To simplify the discussion later, some concepts are briefly explained below. Let $f\left( x \right)$ be a function, the entire real line is denoted by $\mathbf{R}$, $\mathbf{I}$ is a closed interval, and $\mathbf{I}\subseteq \mathbf{R}$. $f\left( x \right)\in {{C}^{n}}\left( \mathbf{I} \right)$ means that the function ${{f}^{\left( n \right)}}\left( x \right)$ is continuous on the closed interval $\mathbf{I}$, that is, the function ${{f}^{\left( n \right)}}\left( x \right)$ is continuous in each of the interior points of the interval $\mathbf{I}$, and is also continuous at the two boundary points of the interval $\mathbf{I}$ respectively, i.e. right continuous and left continuous.

According to theorem\ref{thm2.1}, if ${{R}_{n}}\left[ f \right]=0$, then the function $f\left( x \right)$ satisfies the following equations:
\[{{Z}_{n}}\left[ f \right]-f\left( x \right)=0\]	
Conversely, if the function $f\left( x \right)$ satisfies the equation ${{Z}_{n}}\left[ f \right]-f\left( x \right)=0$, then we have
\[{{R}_{n}}\left[ f \right]=0\]

Therefore, we get the following theorem ( let $y=f\left( x \right)$ ).
\begin{Theorem}\label{thm2.2}
For a given interval $\mathbf{I}$, $a\in \mathbf{I}$, $x\in \mathbf{I}$, the function $y\left( x \right)\in {{C}^{n}}\left( \mathbf{I} \right)$, then the sufficient and necessary condition of function $y\left( x \right)$ satisfies the following linear ordinary differential equation with order $n$.
\begin{equation}\label{eq2.10}
	\sum\limits_{k=1}^{n}{{{\left( -1 \right)}^{k}}\frac{{{\left( x-a \right)}^{k}}}{k!}}{{y}^{(k)}}+y-y\left( a \right)=0
\end{equation}
is
\begin{equation}\label{eq2.11}
	{{R}_{n}}\left( x \right)=\frac{1}{n!}\int\limits_{a}^{x}{{{\left( a-t \right)}^{n}}{{y}^{\left( n+1 \right)}}\left( t \right)dt}=0
\end{equation}
\end{Theorem}

Theorem \ref{thm2.2} points out that for any given function $y\left( x \right)$, if $y\left( x \right)$ satisfies the condition \eqref{eq2.11}, then $y\left( x \right)$ must be a solution to equation \eqref{eq2.10}, which is a linear ordinary differential equation with order $n$.

Let $n\to \infty $, we have the following theorem.
\begin{Theorem}\label{thm2.3}
For a given interval $\mathbf{I}$, $x\in \mathbf{I}$, $a\in \mathbf{I}$, the function $y\left( x \right)\in {{C}^{\infty }}\left( \mathbf{I} \right)$, then the sufficient and necessary condition of function $y\left( x \right)$ satisfies the following linear ordinary differential equation with infinite order
\begin{equation}\label{eq2.12}
	\sum\limits_{k=1}^{\infty }{{{\left( -1 \right)}^{k}}\frac{{{\left( x-a \right)}^{k}}}{k!}}{{y}^{(k)}}+y-y\left( a \right)=0	
\end{equation}
is
\begin{equation}\label{eq2.13}
	\underset{n\to \infty }{\mathop{\lim }}\,{{R}_{n}}\left( x \right)=\underset{n\to \infty }{\mathop{\lim }}\,\frac{1}{n!}\int\limits_{a}^{x}{{{\left( a-t \right)}^{n}}{{y}^{\left( n+1 \right)}}\left( t \right)dt}=0	
\end{equation}
\end{Theorem}

Theorem \ref{thm2.3} points out that for any given function $y\left( x \right)$, if $y\left( x \right)$ satisfies condition \eqref{eq2.13}, then $y\left( x \right)$ must be a solution to equation \eqref{eq2.12}, which is a linear ordinary differential equation with infinite order.

Theorem \ref{thm2.2} and theorem \ref{thm2.3} are thoese of the most important innovations in this paper. From the point of view of application, theorem \ref{thm2.3} is more convenient. 

In some cases, for a given interval $\mathbf{I}$, let $n\to\infty $, if we can get ${{{y}^{\left( n+1 \right)}}\left( x \right)}/{n!}\;\to0$, then condition \eqref{eq2.13} can be satisfied, so, according to theorem \ref{thm2.3}, function  $y\left( x \right)$ is a solution to differential equation \eqref{eq2.12}. For example, for any $x\in \mathbf{R}$, ${{{\lim }_{n\to \infty }}{{\sin }^{\left( n+1 \right)}}\left( \lambda x \right)}/{n!}\;=0$, ${{{\lim }_{n\to \infty }}\exp {{\left( \lambda x \right)}^{\left( n+1 \right)}}}/{n!}\;=0$, where $\lambda $ is a constant, $\lambda \in \mathbf{R}$. So, both $\sin \left( \lambda x \right)$ and ${{e}^{\lambda x}}$ are solutions to differential equation \eqref{eq2.12}.

\begin{Definition}\label{def2.1}
In a given interval $\mathbf{I}$, the set of all solutions satisfying equation \eqref{eq2.12} is called the set of special functions induced by the equation \eqref{eq2.12}, and denoted by ${{Z}^{\infty }}\left( \mathbf{I} \right)$.
\end{Definition}

\begin{Definition}\label{def2.2}
In a given interval $\mathbf{I}$, the set of all solutions satisfying equation \eqref{eq2.10} is called the set of special functions induced by the equation \eqref{eq2.10}, and denoted by ${{Z}^{n}}\left( \mathbf{I} \right)$.
\end{Definition}

Suppose $p\left( x \right)$ is a polynomial, $p\left( x \right)$ obviously satisfies the conditions in theorem \ref{thm2.3}, so $p\left( x \right)$ is a solution to equation \eqref{eq2.12}, i.e. $p\left( x \right)\in {{Z}^{\infty }}\left( \mathbf{I} \right)$. It is easy to prove that both $\sin \left( x \right)$ and ${{e}^{x}}$ satisfy the conditions in theorem \ref{thm2.3}, so $\sin \left( x \right)$ and ${{e}^{x}}$ are also solutions to equation \eqref{eq2.12}, i.e. $\sin \left( x \right)\in {{Z}^{\infty }}\left( \mathbf{I} \right)$, ${{e}^{x}}\in {{Z}^{\infty }}\left( \mathbf{I} \right)$.

Obviously, the functions in set ${{Z}^{\infty }}\left( \mathbf{I} \right)$ are far more than those mentioned above. But how many functions do we have to satisfy equation \eqref{eq2.12}? What role do they play in the classification of functions and in solving differential equations? Those are still to be further studied in the future. 

\section{General theory for solving linear ordinary differential equations}\label{sec3}

This section is the core of this paper. A general method for solving linear ordinary differential equations is studied in this section. Different from the known methods, the new method proposed in this paper uses theorem \ref{thm2.3} to reduce the order by recursive calculations, i.e. we can reduce a given linear ordinary differential equation of order $k$ to the linear ordinary differential equation of order $k-1$. By repeating this recursive process, we can finally get the solution of the given differential equation. 

The core problem of this paper is to find solutions to the following $k$th order linear ordinary differential equation.	
\begin{equation}\label{eq3.1}
{{y}^{\left( k \right)}}={{p}_{k1}}{{y}^{\left( k-1 \right)}}+{{p}_{k2}}{{y}^{\left( k-2 \right)}}+\cdots +{{p}_{k\left( k-1 \right)}}{y}'+{{p}_{kk}}y+{{p}_{k\left( k+1 \right)}}	
\end{equation}
where ${{p}_{k1}}$, ${{p}_{k2}}$, $\cdots$ are all known functions of $x$. Of course, ${{p}_{k1}}$, ${{p}_{k2}}$, $\cdots$, ${{p}_{k\left( k+1 \right)}}$ can also be constants. According to the method given in this paper, a unified solution process will be obtained whether ${{p}_{k1}}$, ${{p}_{k2}}$, $\cdots$, ${{p}_{k\left( k+1 \right)}}$ are the known functions of $x$ or the given constants.

In many practical problems, equation \eqref{eq3.1} needs to satisfy some conditions.

There are many different types of conditions. This paper assumes that the solution to equation \eqref{eq3.1} needs to satisfy the following conditions.
\begin{equation}\label{eq3.2}
\left\{ \begin{aligned}
	& x={{x}_{1}},\ y=y\left( {{x}_{1}} \right) \\ 
	& x={{x}_{2}},\ y=y\left( {{x}_{2}} \right) \\ 
	& \cdots \cdots  \\ 
	& x={{x}_{k}},\ y=y\left( {{x}_{k}} \right) \\ 
\end{aligned} \right.
\end{equation}

I would like to say a few words about the conditions. Many literatures refer to the boundary conditions and initial conditions of differential equations (see \cite{jack_k._hale_ordinary_1980},  \cite{william_a._adkins_ordinary_2010},  \cite{shair_ahmad_textbook_2014} and \cite{mircea_v._soare_ordinary_2007}, etc.). Because the methods in this paper is different from the traditional methods, the boundary condition or initial condition is not used in this paper. In this paper, the set of above conditions \eqref{eq3.2} is called \emph{function value condition} or \emph{condition of function value}.

In general, the $k$th order ordinary differential equation \eqref{eq3.1} satisfying the function value condition \eqref{eq3.2} is denoted as
\begin{equation}\label{eq3.3}
\left\{ \begin{aligned}
	& {{y}^{\left( k \right)}}={{p}_{k1}}{{y}^{\left( k-1 \right)}}+{{p}_{k2}}{{y}^{\left( k-2 \right)}}+\cdots +{{p}_{k\left( k-1 \right)}}{y}'+{{p}_{kk}}y+{{p}_{k\left( k+1 \right)}} \\ 
	& y\left| _{x={{x}_{1}}} \right.=y\left( {{x}_{1}} \right) \\ 
	& y\left| _{x={{x}_{2}}} \right.=y\left( {{x}_{2}} \right) \\ 
	& \cdots \cdots  \\ 
	& y\left| _{x={{x}_{k}}} \right.=y\left( {{x}_{k}} \right) \\ 
\end{aligned} \right.
\end{equation}
where if $i\ne j$, then ${{x}_{i}}\ne {{x}_{j}}$. Equation \eqref{eq3.3} is called \emph{function value problem of order $k$}.

In this paper, we hope to find a general method for solving the function value problem of order $k$ \eqref{eq3.3}.

Find the derivative on both sides of equation \eqref{eq3.1}, we have
\begin{equation}\label{eq3.4}
\begin{aligned}
	{{y}^{\left( k+1 \right)}}& ={{p}_{k1}}\left( {{p}_{k1}}{{y}^{\left( k-1 \right)}}+{{p}_{k2}}{{y}^{\left( k-2 \right)}}+{{p}_{k3}}{{y}^{\left( k-3 \right)}}+\cdots +{{p}_{k\left( k-1 \right)}}{y}'+{{p}_{kk}}y+{{p}_{k\left( k+1 \right)}} \right) \\ 
	& +\left( {{{{p}'}}_{k1}}+{{p}_{k2}} \right){{y}^{\left( k-1 \right)}}+\left( {{{{p}'}}_{k2}}+{{p}_{k3}} \right){{y}^{\left( k-2 \right)}} \\ 
	& +\left( {{{{p}'}}_{k3}}+{{p}_{k4}} \right){{y}^{\left( k-3 \right)}}+\cdots +\left( {{{{p}'}}_{k\left( k-1 \right)}}+{{p}_{kk}} \right){y}'+{{{{p}'}}_{kk}}y+{{{{p}'}}_{k\left( k+1 \right)}} \\ 
	& =\left( {{{{p}'}}_{k1}}+{{p}_{k1}}{{p}_{k1}}+{{p}_{k2}} \right){{y}^{\left( k-1 \right)}}+\left( {{{{p}'}}_{k2}}+{{p}_{k1}}{{p}_{k2}}+{{p}_{k3}} \right){{y}^{\left( k-2 \right)}} \\ 
	& +\left( {{{{p}'}}_{k3}}+{{p}_{k1}}{{p}_{k3}}+{{p}_{k4}} \right){{y}^{\left( k-3 \right)}}+\cdots +\left( {{{{p}'}}_{k\left( k-1 \right)}}+{{p}_{k1}}{{p}_{k\left( k-1 \right)}}+{{p}_{kk}} \right){y}' \\ 
	& +\left( {{{{p}'}}_{kk}}+{{p}_{k1}}{{p}_{kk}} \right)y+\left( {{{{p}'}}_{k\left( k+1 \right)}}+{{p}_{k1}}{{p}_{k\left( k+1 \right)}} \right)  
\end{aligned}
\end{equation}

The above formula can be written as: 
\begin{equation}\label{eq3.5}
{{y}^{\left( k+1 \right)}}={{p}_{\left( k+1 \right)1}}{{y}^{\left( k-1 \right)}}+{{p}_{\left( k+1 \right)2}}{{y}^{\left( k-2 \right)}}+{{p}_{\left( k+1 \right)\left( k-1 \right)}}{y}'+{{p}_{\left( k+1 \right)k}}y+{{p}_{\left( k+1 \right)\left( k+1 \right)}}
\end{equation}
where
\begin{equation}\label{eq3.6}
\left\{ \begin{aligned}
	& {{p}_{\left( k+1 \right)1}}={{{{p}'}}_{k1}}+{{p}_{k1}}{{p}_{k1}}+{{p}_{k2}} \\ 
	& {{p}_{\left( k+1 \right)2}}={{{{p}'}}_{k2}}+{{p}_{k1}}{{p}_{k2}}+{{p}_{k3}} \\ 
	& \cdots \cdots \\ 
	& {{p}_{\left( k+1 \right)\left( k-1 \right)}}={{{{p}'}}_{k\left( k-1 \right)}}+{{p}_{k1}}{{p}_{k\left( k-1 \right)}}+{{p}_{kk}} \\ 
	& {{p}_{\left( k+1 \right)k}}={{{{p}'}}_{kk}}+{{p}_{k1}}{{p}_{kk}} \\ 
	& {{p}_{\left( k+1 \right)\left( k+1 \right)}}={{{{p}'}}_{k\left( k+1 \right)}}+{{p}_{k1}}{{p}_{k\left( k+1 \right)}} \\ 
\end{aligned} \right.
\end{equation}

Derivation on both sides of formula \eqref{eq3.5}, we get
\begin{equation}\label{eq3.7}
{{y}^{\left( k+2 \right)}}={{p}_{\left( k+2 \right)1}}{{y}^{\left( k-1 \right)}}+{{p}_{\left( k+2 \right)2}}{{y}^{\left( k-2 \right)}}+{{p}_{\left( k+2 \right)\left( k-1 \right)}}{y}'+{{p}_{\left( k+2 \right)k}}y+{{p}_{\left( k+2 \right)\left( k+1 \right)}}
\end{equation}
where
\begin{equation}\label{eq3.8}
\left\{ \begin{aligned}
	& {{p}_{\left( k+2 \right)1}}={{{{p}'}}_{\left( k+1 \right)1}}+{{p}_{\left( k+1 \right)1}}{{p}_{k1}}+{{p}_{\left( k+1 \right)2}} \\ 
	& {{p}_{\left( k+2 \right)2}}={{{{p}'}}_{\left( k+1 \right)2}}+{{p}_{\left( k+1 \right)1}}{{p}_{k2}}+{{p}_{\left( k+1 \right)3}} \\ 
	& \cdots \cdots  \\ 
	& {{p}_{\left( k+2 \right)\left( k-1 \right)}}={{{{p}'}}_{\left( k+1 \right)\left( k-1 \right)}}+{{p}_{\left( k+1 \right)1}}{{p}_{\left( k+1 \right)\left( k-1 \right)}}+{{p}_{\left( k+1 \right)k}} \\ 
	& {{p}_{\left( k+2 \right)k}}={{{{p}'}}_{\left( k+1 \right)k}}+{{p}_{\left( k+1 \right)1}}{{p}_{\left( k+1 \right)k}} \\ 
	& {{p}_{\left( k+2 \right)\left( k+1 \right)}}={{{{p}'}}_{\left( k+1 \right)\left( k+1 \right)}}+{{p}_{\left( k+1 \right)1}}{{p}_{\left( k+1 \right)\left( k+1 \right)}} \\ 
\end{aligned} \right.
\end{equation}

Go ahead with the above method and get the following result:
\begin{equation}\label{eq3.9}
{{y}^{\left( k+m \right)}}={{p}_{\left( k+m \right)1}}{{y}^{\left( k-1 \right)}}+{{p}_{\left( k+m \right)2}}{{y}^{\left( k-2 \right)}}+{{p}_{\left( k+m \right)\left( k-1 \right)}}{y}'+{{p}_{\left( k+m \right)k}}y+{{p}_{\left( k+m \right)\left( k+1 \right)}}
\end{equation}
where $m=1,2,\cdots $, and
\begin{equation}\label{eq3.10}
\left\{ \begin{aligned}
	& {{p}_{\left( k+m \right)1}}={{{{p}'}}_{\left( k+m-1 \right)1}}+{{p}_{\left( k+m-1 \right)1}}{{p}_{k1}}+{{p}_{\left( k+m-1 \right)2}} \\ 
	& {{p}_{\left( k+m \right)2}}={{{{p}'}}_{\left( k+m-1 \right)2}}+{{p}_{\left( k+m-1 \right)1}}{{p}_{k2}}+{{p}_{\left( k+m-1 \right)3}} \\ 
	& \cdots \cdots \\ 
	& {{p}_{\left( k+m \right)\left( k-1 \right)}}={{{{p}'}}_{\left( k+m-1 \right)\left( k-1 \right)}}+{{p}_{\left( k+m-1 \right)1}}{{p}_{\left( k+m-1 \right)\left( k-1 \right)}}+{{p}_{\left( k+m-1 \right)k}} \\ 
	& {{p}_{\left( k+m \right)k}}={{{{p}'}}_{\left( k+m-1 \right)k}}+{{p}_{\left( k+m-1 \right)1}}{{p}_{\left( k+m-1 \right)k}} \\ 
	& {{p}_{\left( k+m \right)\left( k+1 \right)}}={{{{p}'}}_{\left( k+m-1 \right)\left( k+1 \right)}}+{{p}_{\left( k+m-1 \right)1}}{{p}_{\left( k+m-1 \right)\left( k+1 \right)}} \\ 
\end{aligned} \right.
\end{equation}

Let
\begin{equation}\label{eq3.11}
n=k+m
\end{equation}
then the derivative formula \eqref{eq3.9} can be written as follows:
\begin{equation}\label{eq3.12}
{{y}^{\left( n \right)}}={{p}_{n1}}{{y}^{\left( k-1 \right)}}+{{p}_{n2}}{{y}^{\left( k-2 \right)}}+{{p}_{n\left( k-1 \right)}}{y}'+{{p}_{nk}}y+{{p}_{n\left( k+1 \right)}}
\end{equation}

The derivative of the $n$th order can be written as a linear combination of ${{y}^{\left( k-1 \right)}}$, ${{y}^{\left( k-2 \right)}}$, $\cdots$, ${y}'$ and $y$. Its coefficients ${{p}_{n1}}$, ${{p}_{n2}}$, $\cdots$, ${{p}_{n\left( k+1 \right)}}$ are functions of $x$ and  can be obtained by the following recursive formulas:
\begin{equation}\label{eq3.13}
\left\{ \begin{aligned}
	& {{p}_{n1}}={{{{p}'}}_{\left( n-1 \right)1}}+{{p}_{\left( n-1 \right)1}}{{p}_{k1}}+{{p}_{\left( n-1 \right)2}} \\ 
	& {{p}_{n2}}={{{{p}'}}_{\left( n-1 \right)2}}+{{p}_{\left( n-1 \right)1}}{{p}_{k2}}+{{p}_{\left( n-1 \right)3}} \\ 
	& \cdots \cdots \\ 
	& {{p}_{n\left( k-1 \right)}}={{{{p}'}}_{\left( n-1 \right)\left( k-1 \right)}}+{{p}_{\left( n-1 \right)1}}{{p}_{\left( n-1 \right)\left( k-1 \right)}}+{{p}_{\left( n-1 \right)k}} \\ 
	& {{p}_{nk}}={{{{p}'}}_{\left( n-1 \right)k}}+{{p}_{\left( n-1 \right)1}}{{p}_{\left( n-1 \right)k}} \\ 
	& {{p}_{n\left( k+1 \right)}}={{{{p}'}}_{\left( n-1 \right)\left( k+1 \right)}}+{{p}_{\left( n-1 \right)1}}{{p}_{\left( n-1 \right)\left( k+1 \right)}} \\ 
\end{aligned} \right.
\end{equation}
where $n=k+1,k+2,\cdots $. For linear ordinary differential equations of order $k$, ${{p}_{k1}}$, ${{p}_{k2}}$, $\cdots$, ${{p}_{k\left( k+1 \right)}}$ are known quantities. 

By substituting the last expression $y\left| _{x={{x}_{k}}} \right.=y\left( {{x}_{k}} \right)$ of the function value condition in \eqref{eq3.3} into formula \eqref{eq2.12}, we have
\begin{equation}\label{eq3.14}
\sum\limits_{n=1}^{k-1}{{{\left( -1 \right)}^{n}}\frac{{{\left( x-{{x}_{k}} \right)}^{n}}}{n!}}{{y}^{(n)}}+\sum\limits_{n=k}^{\infty }{{{\left( -1 \right)}^{n}}\frac{{{\left( x-{{x}_{k}} \right)}^{n}}}{n!}}{{y}^{(n)}}+y-y\left( {{x}_{k}} \right)=0
\end{equation}

Substituting \eqref{eq3.1} into \eqref{eq3.14}, we get
\begin{equation}\label{eq3.15}
\begin{aligned}
	& \sum\limits_{n=k}^{\infty }{{{\left( -1 \right)}^{n}}\frac{{{\left( x-{{x}_{k}} \right)}^{n}}}{n!}}\left( {{p}_{n1}}{{y}^{\left( k-1 \right)}}+{{p}_{n2}}{{y}^{\left( k-2 \right)}}+\cdots +{{p}_{n\left( k-1 \right)}}{y}'+{{p}_{nk}}y+{{p}_{n\left( k+1 \right)}} \right) \\ 
	& +\sum\limits_{n=1}^{k-1}{{{\left( -1 \right)}^{n}}\frac{{{\left( x-{{x}_{k}} \right)}^{n}}}{n!}}{{y}^{(n)}}+y-y\left( {{x}_{k}} \right)=0 \\ 
\end{aligned}
\end{equation}
i.e.
\begin{equation}\label{eq3.16}
\begin{aligned}
	& \left( \sum\limits_{n=k}^{\infty }{{{\left( -1 \right)}^{n}}\frac{{{\left( x-{{x}_{k}} \right)}^{n}}}{n!}}{{p}_{n1}}+{{\left( -1 \right)}^{k-1}}\frac{{{\left( x-{{x}_{k}} \right)}^{k-1}}}{\left( k-1 \right)!} \right){{y}^{\left( k-1 \right)}} \\ 
	& +\left( \sum\limits_{n=k}^{\infty }{{{\left( -1 \right)}^{n}}\frac{{{\left( x-{{x}_{k}} \right)}^{n}}}{n!}}{{p}_{n2}}+{{\left( -1 \right)}^{k-2}}\frac{{{\left( x-{{x}_{k}} \right)}^{k-2}}}{\left( k-2 \right)!} \right){{y}^{\left( k-2 \right)}} \\ 
	& +\cdots \cdots  \\ 
	& +\left( \sum\limits_{n=k}^{\infty }{{{\left( -1 \right)}^{n}}\frac{{{\left( x-{{x}_{k}} \right)}^{n}}}{n!}}{{p}_{n\left( k-1 \right)}}+{{\left( -1 \right)}^{1}}\frac{\left( x-{{x}_{k}} \right)}{1!} \right){y}' \\ 
	& +\left( \sum\limits_{n=k}^{\infty }{{{\left( -1 \right)}^{n}}\frac{{{\left( x-{{x}_{k}} \right)}^{n}}}{n!}}{{p}_{nk}}+{{\left( -1 \right)}^{0}}\frac{{{\left( x-{{x}_{k}} \right)}^{0}}}{0!} \right)y \\ 
	& +\left( \sum\limits_{n=k}^{\infty }{{{\left( -1 \right)}^{n}}\frac{{{\left( x-{{x}_{k}} \right)}^{n}}}{n!}}{{p}_{n\left( k+1 \right)}}-y\left( {{x}_{k}} \right) \right)=0 \\ 
\end{aligned}
\end{equation}

Therefore, the linear ordinary differential equation of  order $k$ can be reduced to the above linear ordinary differential equation of order $k-1$.

The linear ordinary differential equations of order $k-1$ \eqref{eq3.16} can be written in the following form:
\begin{equation}\label{eq3.17}
{{Q}_{\left( k-1 \right)1}}{{y}^{\left( k-1 \right)}}+{{Q}_{\left( k-1 \right)2}}{{y}^{\left( k-2 \right)}}+\cdots +{{Q}_{\left( k-1 \right)\left( k-1 \right)}}{y}'+{{Q}_{\left( k-1 \right)k}}y+{{Q}_{\left( k-1 \right)\left( k+1 \right)}}=0
\end{equation}

Equation \eqref{eq3.17} has automatically contained an function value condition $y\left| _{x={{x}_{k}}} \right.=y\left( {{x}_{k}} \right)$, where
\begin{equation}\label{eq3.18}
{{Q}_{\left( k-1 \right)1}}=\sum\limits_{n=k}^{\infty }{{{\left( -1 \right)}^{n}}\frac{{{\left( x-{{x}_{k}} \right)}^{n}}}{n!}}{{p}_{n1}}+{{\left( -1 \right)}^{k-1}}\frac{{{\left( x-{{x}_{k}} \right)}^{k-1}}}{\left( k-1 \right)!}
\end{equation}
\begin{equation}\label{eq3.19}
{{Q}_{\left( k-1 \right)2}}=\sum\limits_{n=k}^{\infty }{{{\left( -1 \right)}^{n}}\frac{{{\left( x-{{x}_{k}} \right)}^{n}}}{n!}}{{p}_{n2}}+{{\left( -1 \right)}^{k-2}}\frac{{{\left( x-{{x}_{k}} \right)}^{k-2}}}{\left( k-2 \right)!}
\end{equation}
\begin{equation*}
\cdots \cdots
\end{equation*}
\begin{equation}\label{eq3.20}
{{Q}_{\left( k-1 \right)\left( k-1 \right)}}=\sum\limits_{n=k}^{\infty }{{{\left( -1 \right)}^{n}}\frac{{{\left( x-{{x}_{k}} \right)}^{n}}}{n!}}{{p}_{n\left( k-1 \right)}}+{{\left( -1 \right)}^{1}}\frac{\left( x-{{x}_{k}} \right)}{1!}
\end{equation}
\begin{equation}\label{eq3.21}
{{Q}_{\left( k-1 \right)k}}=\sum\limits_{n=k}^{\infty }{{{\left( -1 \right)}^{n}}\frac{{{\left( x-{{x}_{k}} \right)}^{n}}}{n!}}{{p}_{nk}}+{{\left( -1 \right)}^{0}}\frac{{{\left( x-{{x}_{k}} \right)}^{0}}}{0!}
\end{equation}
\begin{equation}\label{eq3.22}
{{Q}_{\left( k-1 \right)\left( k+1 \right)}}=\sum\limits_{n=k}^{\infty }{{{\left( -1 \right)}^{n}}\frac{{{\left( x-{{x}_{k}} \right)}^{n}}}{n!}}{{p}_{n\left( k+1 \right)}}-y\left( {{x}_{k}} \right)
\end{equation}

Thus, when ${{Q}_{\left( k-1 \right)1}}\ne 0$, the $\left( k-1 \right)$th order linear differential equation \eqref{eq3.17} can be written in the following form: 
\begin{equation}\label{eq3.23}
{{y}^{\left( k-1 \right)}}=-\frac{{{Q}_{\left( k-1 \right)2}}}{{{Q}_{\left( k-1 \right)1}}}{{y}^{\left( k-2 \right)}}-\frac{{{Q}_{\left( k-1 \right)3}}}{{{Q}_{\left( k-1 \right)1}}}{{y}^{\left( k-3 \right)}}-\cdots -\frac{{{Q}_{\left( k-1 \right)\left( k-1 \right)}}}{{{Q}_{\left( k-1 \right)1}}}{y}'-\frac{{{Q}_{\left( k-1 \right)k}}}{{{Q}_{\left( k-1 \right)1}}}y-\frac{{{Q}_{\left( k-1 \right)\left( k+1 \right)}}}{{{Q}_{\left( k-1 \right)1}}}
\end{equation}

Let
\begin{equation}\label{eq3.24}
\left\{ \begin{aligned}
	& {{p}_{\left( k-1 \right)1}}=-\frac{{{Q}_{\left( k-1 \right)2}}}{{{Q}_{\left( k-1 \right)1}}} \\ 
	& {{p}_{\left( k-1 \right)2}}=-\frac{{{Q}_{\left( k-1 \right)3}}}{{{Q}_{\left( k-1 \right)1}}} \\ 
	& \cdots \cdots  \\ 
	& {{p}_{\left( k-1 \right)\left( k-2 \right)}}=-\frac{{{Q}_{\left( k-1 \right)\left( k-1 \right)}}}{{{Q}_{\left( k-1 \right)1}}} \\ 
	& {{p}_{\left( k-1 \right)\left( k-1 \right)}}=-\frac{{{Q}_{\left( k-1 \right)k}}}{{{Q}_{\left( k-1 \right)1}}} \\ 
	& {{p}_{\left( k-1 \right)k}}=-\frac{{{Q}_{\left( k-1 \right)\left( k+1 \right)}}}{{{Q}_{\left( k-1 \right)1}}} \\ 
\end{aligned} \right.
\end{equation}
then equation \eqref{eq3.23} can be written in the following:
\begin{equation}\label{eq3.25}
{{y}^{\left( k-1 \right)}}={{p}_{\left( k-1 \right)1}}{{y}^{\left( k-2 \right)}}+{{p}_{\left( k-1 \right)2}}{{y}^{\left( k-3 \right)}}+\cdots +{{p}_{\left( k-1 \right)\left( k-2 \right)}}{y}'+{{p}_{\left( k-1 \right)\left( k-1 \right)}}y+{{p}_{\left( k-1 \right)k}}
\end{equation}

By comparing the above equation with equation \eqref{eq3.1}, we know that equation \eqref{eq3.25} is a $\left( k-1 \right)$th order linear ordinary differential equation after the reduction of order by the $k$th order linear ordinary differential equation \eqref{eq3.1}.

Since the equation \eqref{eq3.25} has automatically contained an function value condition $y\left| _{x={{x}_{k}}} \right.=y\left( {{x}_{k}} \right)$, the original $k$ function value conditions are reduced by 1, and become the $k-1$ function value conditions, or the original $k$ function value conditions are reduced to $k-1$ function value conditions. Therefore, by reducing order, the original function value problem of order $k$ \eqref{eq3.3} becomes the following function value problem order $k-1$:
\begin{equation}\label{eq3.26}
\left\{ \begin{aligned}
	& {{y}^{\left( k-1 \right)}}={{p}_{\left( k-1 \right)1}}{{y}^{\left( k-2 \right)}}+{{p}_{\left( k-1 \right)2}}{{y}^{\left( k-3 \right)}}+\cdots +{{p}_{\left( k-1 \right)\left( k-2 \right)}}{y}'+{{p}_{\left( k-1 \right)\left( k-1 \right)}}y+{{p}_{\left( k-1 \right)k}} \\ 
	& y\left| _{x={{x}_{1}}} \right.=y\left( {{x}_{1}} \right) \\ 
	& y\left| _{x={{x}_{2}}} \right.=y\left( {{x}_{2}} \right) \\ 
	& \cdots \cdots  \\ 
	& y\left| _{x={{x}_{k-1}}} \right.=y\left( {{x}_{k-1}} \right) \\ 
\end{aligned} \right.
\end{equation}

With the same method, the$\left( k-1 \right)$th order linear ordinary differential equation \eqref{eq3.25} can be reduced to the $\left( k-2 \right)$th order linear ordinary differential equation, but we need to substitute the function value condition $y\left| _{x={{x}_{k-1}}} \right.=y\left( {{x}_{k-1}} \right)$ into the following formula (see \eqref{eq2.12}):
\begin{equation}\label{eq3.27}
\sum\limits_{n=1}^{k-2}{{{\left( -1 \right)}^{n}}\frac{{{\left( x-{{x}_{k-1}} \right)}^{n}}}{n!}}{{y}^{(n)}}+\sum\limits_{n=k-1}^{\infty }{{{\left( -1 \right)}^{n}}\frac{{{\left( x-{{x}_{k-1}} \right)}^{n}}}{n!}}{{y}^{(n)}}+y-y\left( {{x}_{k-1}} \right)=0
\end{equation}

By deriving the formula similar to the previous one, we can get the following result:
\begin{equation}\label{eq3.28}
\left\{ \begin{aligned}
	& {{y}^{\left( k-2 \right)}}={{p}_{\left( k-2 \right)1}}{{y}^{\left( k-3 \right)}}+{{p}_{\left( k-2 \right)2}}{{y}^{\left( k-4 \right)}}+\cdots +{{p}_{\left( k-2 \right)\left( k-3 \right)}}{y}'+{{p}_{\left( k-2 \right)\left( k-2 \right)}}y+{{p}_{\left( k-2 \right)\left( k-1 \right)}} \\ 
	& y\left| _{x={{x}_{1}}} \right.=y\left( {{x}_{1}} \right) \\ 
	& y\left| _{x={{x}_{2}}} \right.=y\left( {{x}_{2}} \right) \\ 
	& \cdots \cdots  \\ 
	& y\left| _{x={{x}_{k-2}}} \right.=y\left( {{x}_{k-2}} \right) \\ 
\end{aligned} \right.
\end{equation}

The two function value conditions $y\left| _{x={{x}_{k-1}}} \right.=y\left( {{x}_{k-1}} \right)$ and $y\left| _{x={{x}_{k-2}}} \right.=y\left( {{x}_{k-2}} \right)$ have been included in equation \eqref{eq3.28}, where
\begin{equation}\label{eq3.29}
\left\{ \begin{aligned}
	& {{p}_{\left( k-2 \right)1}}=-\frac{{{Q}_{\left( k-2 \right)2}}}{{{Q}_{\left( k-2 \right)1}}} \\ 
	& {{p}_{\left( k-2 \right)2}}=-\frac{{{Q}_{\left( k-2 \right)3}}}{{{Q}_{\left( k-2 \right)1}}} \\ 
	& \cdots \cdots  \\ 
	& {{p}_{\left( k-2 \right)\left( k-3 \right)}}=-\frac{{{Q}_{\left( k-2 \right)\left( k-2 \right)}}}{{{Q}_{\left( k-2 \right)1}}} \\ 
	& {{p}_{\left( k-2 \right)\left( k-2 \right)}}=-\frac{{{Q}_{\left( k-2 \right)\left( k-1 \right)}}}{{{Q}_{\left( k-2 \right)1}}} \\ 
	& {{p}_{\left( k-2 \right)\left( k-1 \right)}}=\frac{{{Q}_{\left( k-2 \right)k}}}{{{Q}_{\left( k-2 \right)1}}} \\ 
\end{aligned} \right.
\end{equation}
where ${{Q}_{\left( k-2 \right)1}}\ne 0$, and
\begin{equation}\label{eq3.30}
{{Q}_{\left( k-2 \right)1}}=\sum\limits_{n=k-1}^{\infty }{{{\left( -1 \right)}^{n}}\frac{{{\left( x-{{x}_{k-1}} \right)}^{n}}}{n!}}{{p}_{n1}}+{{\left( -1 \right)}^{k-2}}\frac{{{\left( x-{{x}_{k-1}} \right)}^{k-2}}}{\left( k-2 \right)!}
\end{equation}
\begin{equation}\label{eq3.31}
{{Q}_{\left( k-2 \right)2}}=\sum\limits_{n=k-1}^{\infty }{{{\left( -1 \right)}^{n}}\frac{{{\left( x-{{x}_{k-1}} \right)}^{n}}}{n!}}{{p}_{n2}}+{{\left( -1 \right)}^{k-3}}\frac{{{\left( x-{{x}_{k-1}} \right)}^{k-3}}}{\left( k-3 \right)!}
\end{equation}
\[\cdots \cdots\]
\begin{equation}\label{eq3.32}
{{Q}_{\left( k-2 \right)\left( k-2 \right)}}=\sum\limits_{n=k-1}^{\infty }{{{\left( -1 \right)}^{n}}\frac{{{\left( x-{{x}_{k-1}} \right)}^{n}}}{n!}}{{p}_{n\left( k-2 \right)}}+{{\left( -1 \right)}^{1}}\frac{\left( x-{{x}_{k-1}} \right)}{1!}
\end{equation}
\begin{equation}\label{eq3.33}
{{Q}_{\left( k-2 \right)\left( k-1 \right)}}=\sum\limits_{n=k-1}^{\infty }{{{\left( -1 \right)}^{n}}\frac{{{\left( x-{{x}_{k-1}} \right)}^{n}}}{n!}}{{p}_{n\left( k-1 \right)}}+{{\left( -1 \right)}^{0}}\frac{{{\left( x-{{x}_{k-1}} \right)}^{0}}}{0!}
\end{equation}
\begin{equation}\label{eq3.34}
{{Q}_{\left( k-2 \right)k}}=\sum\limits_{n=k-1}^{\infty }{{{\left( -1 \right)}^{n}}\frac{{{\left( x-{{x}_{k-1}} \right)}^{n}}}{n!}}{{p}_{nk}}-y\left( {{x}_{k-1}} \right)
\end{equation}
\begin{equation}\label{eq3.35}
\left\{ \begin{aligned}
	& {{p}_{n1}}={{{{p}'}}_{\left( n-1 \right)1}}+{{p}_{\left( n-1 \right)1}}{{p}_{\left( k-1 \right)1}}+{{p}_{\left( n-1 \right)2}} \\ 
	& {{p}_{n2}}={{{{p}'}}_{\left( n-1 \right)2}}+{{p}_{\left( n-1 \right)1}}{{p}_{\left( k-1 \right)2}}+{{p}_{\left( n-1 \right)3}} \\ 
	& \cdots \cdots \\ 
	& {{p}_{n\left( k-2 \right)}}={{{{p}'}}_{\left( n-1 \right)\left( k-2 \right)}}+{{p}_{\left( n-1 \right)1}}{{p}_{\left( k-1 \right)\left( k-2 \right)}}+{{p}_{\left( n-1 \right)\left( k-1 \right)}} \\ 
	& {{p}_{n\left( k-1 \right)}}={{{{p}'}}_{\left( n-1 \right)\left( k-1 \right)}}+{{p}_{\left( n-1 \right)1}}{{p}_{\left( k-1 \right)\left( k-1 \right)}} \\ 
	& {{p}_{nk}}={{{{p}'}}_{\left( n-1 \right)k}}+{{p}_{\left( n-1 \right)1}}{{p}_{\left( k-1 \right)k}} \\ 
\end{aligned} \right.
\end{equation}
where $n=k$, $k+1$, $\cdots$. For $\left( k-1 \right)$th order linear differential equations, ${{p}_{\left( k-1 \right)1}}$, ${{p}_{\left( k-1 \right)2}}$, $\cdots$, ${{p}_{\left( k-1 \right)k}}$ are all known quantities, which is calculated by formula \eqref{eq3.24}.

After a series of reduction of order in the same way, we can find the second order linear ordinary differential equations, the first order linear ordinary differential equations and the 0th order linear ordinary differential equations in the following.

The second order linear differential equation:
\begin{equation}\label{eq3.36}
\left\{ \begin{aligned}
	& {y}''={{p}_{21}}{y}'+{{p}_{22}}y+{{p}_{23}} \\ 
	& y\left| _{x={{x}_{1}}} \right.=y\left( {{x}_{1}} \right) \\ 
	& y\left| _{x={{x}_{2}}} \right.=y\left( {{x}_{2}} \right) \\ 
\end{aligned} \right.
\end{equation}

The $k-2$ function value conditions $y\left| _{x={{x}_{3}}} \right.=y\left( {{x}_{3}} \right)$, $\cdots$, $y\left| _{x={{x}_{k}}} \right.=y\left( {{x}_{k}} \right)$ have been included in equation \eqref{eq3.36}

The first order linear differential equation:
\begin{equation}\label{eq3.37}
\left\{ \begin{aligned}
	& {y}'={{p}_{11}}y+{{p}_{12}} \\ 
	& y\left| _{x={{x}_{1}}} \right.=y\left( {{x}_{1}} \right) \\ 
\end{aligned} \right.
\end{equation}

The $k-1$ function value conditions $y\left| _{x={{x}_{2}}} \right.=y\left( {{x}_{2}} \right)$, $\cdots$, $y\left| _{x={{x}_{k}}} \right.=y\left( {{x}_{k}} \right)$ have been included in equation \eqref{eq3.37}

The 0th order linear differential equation:
\begin{equation}\label{eq3.38}
y={{p}_{01}}=-\frac{{{Q}_{02}}}{{{Q}_{01}}}=\frac{\sum\limits_{n=1}^{\infty }{{{\left( -1 \right)}^{n}}\frac{{{\left( x-{{x}_{1}} \right)}^{n}}}{n!}{{p}_{n2}}-y\left( {{x}_{1}} \right)}}{\sum\limits_{n=1}^{\infty }{{{\left( -1 \right)}^{n}}\frac{{{\left( x-{{x}_{1}} \right)}^{n}}}{n!}{{p}_{n1}}+1}}
\end{equation}

All $k$ function value conditions $y\left| _{x={{x}_{1}}} \right.=y\left( {{x}_{1}} \right)$, $\cdots$, $y\left| _{x={{x}_{k}}} \right.=y\left( {{x}_{k}} \right)$ have been included in equation \eqref{eq3.38}. Formula \eqref{eq3.38} is the solution to function value problem of order $k$ \eqref{eq3.3}.

The above results can be summed up as the following important theorem.
\begin{Theorem}\label{thm3.1}
Suppose that on a closed interval $\mathbf{I}$, ${{x}_{1}}\in \mathbf{I}$, ${{x}_{2}}\in \mathbf{I}$, $\cdots$, ${{x}_{k}}\in \mathbf{I}$, $a\in \mathbf{I}$, $x\in \mathbf{I}$, and ${{x}_{1}}$, ${{x}_{2}}$, $\cdots$, ${{x}_{k}}$ are different from one another, ${{p}_{k1}}$, ${{p}_{k2}}$, $\cdots$, ${{p}_{k\left( k+1 \right)}}$ are functions of $x$, assume that the solution $y=y\left( x \right)$ of equation \eqref{eq3.3} and ${{p}_{k1}}$, ${{p}_{k2}}$, $\cdots$, ${{p}_{k\left( k+1 \right)}}$ have derivatives of arbitrary order on the interval $\mathbf{I}$, and 
\[\underset{n\to \infty }{\mathop{\lim }}\,{{R}_{n}}\left( x \right)=\underset{n\to \infty }{\mathop{\lim }}\,\frac{1}{n!}\int\limits_{a}^{x}{{{\left( {{x}_{i}}-t \right)}^{n}}{{y}^{\left( n+1 \right)}}\left( t \right)dt}=0\]
where $i=1,2, \cdots,k$, ${{Q}_{\left( k-1 \right)1}}\ne 0$, ${{Q}_{\left( k-2 \right)1}}\ne 0$, $\cdots$, ${{Q}_{01}}\ne 0$, then the function value problem of order $k$ $($see \eqref{eq3.3}$)$ becomes the function value problem of order $k-1$ $($see \eqref{eq3.26}$)$ by \eqref{eq3.13}--\eqref{eq3.16}.
\end{Theorem}

By applying theorem \ref{thm3.1} repeatedly, we can get the solution to the $k$th order linear ordinary differential equation \eqref{eq3.3}.

In the conditions of theorem\ref{thm3.1}, except the condition ${{\lim }_{n\to \infty }}{{R}_{n}}\left( x \right)=0$, the other conditions are easy to be satisfied in many practical problems. It can be assumed that a large number of functions can satisfy the condition ${{\lim }_{n\to \infty }}{{R}_{n}}\left( x \right)=0$, which needs further study. The in-depth discussion of this issue is beyond the scope of this paper. I will give it in one of my follow-up papers. 

It can be seen from the above discussion that the core process of solving linear ordinary differential equations is that the order can be gradually reduced by recursion, and finally, the solution is obtained, and the function value conditions are automatically included in the solution.

\section{Examples}\label{sec4}
The two examples given in this section are examples of the second order linear ordinary differential equations. Because the second order linear ordinary differential equations are widely used (see \cite{jack_k._hale_ordinary_1980}, \cite{william_a._adkins_ordinary_2010},  \cite{shair_ahmad_textbook_2014} and \cite{mircea_v._soare_ordinary_2007}), I have discussed in the earlier version (v1 or v2) of this paper (see \cite{daiyuan_zhang_general_2018}) for the second order linear ordinary differential equations \eqref{eq3.36} in detail. 

Example \ref{exm4.1} is a simple one for sovling the second order linear ordinary differential equation with constant coefficients of known solutions. This example is chosen to verify the correctness of the method proposed in this paper.

Example \ref{exm4.2} is to solve a second order linear ordinary differential equation with variable coefficients. The solving proceses and steps are the same as that of example \ref{exm4.1}. It shows that the method presented in this paper is a general method.  

\subsection{An example for sovling second order linear ordinary differential equation with constant coefficients}\label{subsec4.1}
\begin{Example}{Solve the following equation.}\label{exm4.1}
\begin{equation}\label{eq4.1}
	\left\{ \begin{aligned}
		& {y}''=y \\ 
		& {{\left. y \right|}_{x=b}}=y\left( b \right) \\ 
		& {{\left. y \right|}_{x=a}}=y\left( a \right) \\ 
	\end{aligned} \right.
\end{equation}
\end{Example}

\textbf{Solution}. Obviously, let $k=2$, from equation \eqref{eq3.13}, we have

\begin{equation}\label{eq4.1a}
\left\{ \begin{aligned}
	& {{p}_{n1}}={{{{p}'}}_{\left( n-1 \right)1}}+{{p}_{\left( n-1 \right)1}}{{p}_{21}}+{{p}_{\left( n-1 \right)2}} \\ 
	& {{p}_{n2}}={{{{p}'}}_{\left( n-1 \right)2}}+{{p}_{\left( n-1 \right)1}}{{p}_{22}} \\ 
	& {{p}_{n3}}={{{{p}'}}_{\left( n-1 \right)3}}+{{p}_{\left( n-1 \right)1}}{{p}_{23}} \\ 
\end{aligned} \right.
\end{equation}

From \eqref{eq4.1}, we get	
\begin{equation}\label{eq4.2}
\left\{ \begin{aligned}
	& {{p}_{21}}=0 \\ 
	& {{p}_{22}}=1 \\ 
	& {{p}_{23}}=0 \\ 
\end{aligned} \right.
\end{equation}

Using recursive formulas \eqref{eq4.1a}, we get
\begin{equation}\label{eq4.3}
{{p}_{31}}={{p}_{22}}=1
\end{equation}
\begin{equation}\label{eq4.4}
{{p}_{32}}={{{p}'}_{22}}=0
\end{equation}
\begin{equation}\label{eq4.5}
{{p}_{33}}=0
\end{equation}
\begin{equation}\label{eq4.6}
{{p}_{41}}=2{{{p}'}_{22}}=0
\end{equation}
\begin{equation}\label{eq4.7}
{{p}_{42}}={{{p}''}_{22}}+p_{22}^{2}=1
\end{equation}
\begin{equation}\label{eq4.8}
{{p}_{43}}=0
\end{equation}
\begin{equation}\label{eq4.9}
{{p}_{51}}=3{{{p}''}_{22}}+p_{22}^{2}=1
\end{equation}
\begin{equation}\label{eq4.10}
{{p}_{52}}={{{p}'''}_{22}}+4{{{p}'}_{22}}{{p}_{22}}=0
\end{equation}
\begin{equation}\label{eq4.11}
{{p}_{53}}=0
\end{equation}
\begin{equation}\label{eq4.12}
{{p}_{61}}=4{{{p}'''}_{22}}+6{{{p}'}_{22}}{{p}_{22}}=0
\end{equation}
\begin{equation}\label{eq4.13}
{{p}_{62}}=p_{22}^{\left( 4 \right)}+7{{{p}''}_{22}}{{p}_{22}}+4{{\left( {{{{p}'}}_{22}} \right)}^{2}}+p_{22}^{3}=1
\end{equation}
\begin{equation}\label{eq4.14}
{{p}_{63}}=0
\end{equation}
\begin{equation*}
\cdots \cdots
\end{equation*}

From \eqref{eq3.20}--\eqref{eq3.22}, let $k=2$, ${{x}_{1}}=b$, ${{x}_{2}}=a$, we have
\begin{equation}\label{eq4.15}
\begin{aligned}
	{{Q}_{11}}& =\sum\limits_{n=2}^{\infty }{{{\left( -1 \right)}^{n}}\frac{{{\left( x-a \right)}^{n}}}{n!}}{{p}_{n1}}-\left( x-a \right) \\ 
	& =-\sum\limits_{n=0}^{\infty }{\frac{{{\left( x-a \right)}^{2n+1}}}{\left( 2n+1 \right)!}} \\ 
	& =-\left( \frac{x-a}{1}+\frac{{{\left( x-a \right)}^{3}}}{3!}+\frac{{{\left( x-a \right)}^{5}}}{5!}+\cdots  \right) \\ 
	& =-\sinh \left( x-a \right)  
\end{aligned}
\end{equation}
\begin{equation}\label{eq4.16}
\begin{aligned}
	{{Q}_{12}}& =\sum\limits_{n=2}^{\infty }{{{\left( -1 \right)}^{n}}\frac{{{\left( x-a \right)}^{n}}}{n!}}{{p}_{n2}}+1 \\ 
	& =\sum\limits_{n=0}^{\infty }{\frac{{{\left( x-a \right)}^{2n}}}{\left( 2n \right)!}} \\ 
	& =1+\frac{{{\left( x-a \right)}^{2}}}{2!}+\frac{{{\left( x-a \right)}^{4}}}{4!}+\cdots  \\ 
	& =\cosh \left( x-a \right)  
\end{aligned}
\end{equation}
\begin{equation}\label{eq4.17}
\begin{aligned}
	{{Q}_{13}}& =\sum\limits_{n=2}^{\infty }{{{\left( -1 \right)}^{n}}\frac{{{\left( x-a \right)}^{n}}}{n!}}{{p}_{n3}}-y\left( a \right) =-y\left( a \right)  
\end{aligned}
\end{equation}

After order reduction, we get the first order linear ordinary differential equation in the following: 
\begin{equation}\label{eq4.18}
\left\{ \begin{aligned}
	& {y}'={{p}_{11}}y+{{p}_{12}} \\ 
	& {{\left. y \right|}_{x=b}}=y\left( b \right) \\ 
\end{aligned} \right.
\end{equation}
where, 
\begin{equation}\label{eq4.19}
{{p}_{11}}=-\frac{{{Q}_{12}}\left( x \right)}{{{Q}_{11}}\left( x \right)}=\frac{\cosh\left( x-a \right)}{\sinh\left( x-a \right)}=\coth \left( x-a \right)
\end{equation}
\begin{equation}\label{eq4.20}
{{p}_{12}}=-\frac{{{Q}_{13}}\left( x \right)}{{{Q}_{11}}\left( x \right)}=-\frac{y\left( a \right)}{\sinh\left( x-a \right)}
\end{equation}

From \eqref{eq3.33}--\eqref{eq3.34}, let $k=2$, ${{x}_{1}}=b$, ${{x}_{2}}=a$, we have
\begin{equation}\label{eq4.21}
\begin{aligned}
	{{Q}_{01}}& =\sum\limits_{n=1}^{\infty }{{{\left( -1 \right)}^{n}}\frac{{{\left( x-b \right)}^{n}}}{n!}{{p}_{n1}}+1} \\ 
	& =1-\left( x-b \right){{p}_{11}}+\frac{{{\left( x-b \right)}^{2}}}{2!}{{p}_{21}}-\frac{{{\left( x-b \right)}^{3}}}{3!}{{p}_{31}}+\frac{{{\left( x-b \right)}^{4}}}{4!}{{p}_{41}}-\cdots  \\ 
	& =1-\left( x-b \right)\coth \left( x-a \right)+\frac{{{\left( x-b \right)}^{2}}}{2!}-\frac{{{\left( x-b \right)}^{3}}}{3!}\coth \left( x-a \right)+\frac{{{\left( x-b \right)}^{4}}}{4!}-\cdots  \\ 
	& =1+\frac{{{\left( x-b \right)}^{2}}}{2!}+\frac{{{\left( x-b \right)}^{4}}}{4!}+\cdots -\coth \left( x-a \right)\left( \left( x-b \right)+\frac{{{\left( x-b \right)}^{3}}}{3!}+\cdots  \right) \\ 
	& =\cosh \left( x-b \right)-\coth \left( x-a \right)\sinh \left( x-b \right)  
\end{aligned}
\end{equation}
\begin{equation}\label{eq4.22}
\begin{aligned}
	{{Q}_{02}}& =\sum\limits_{n=1}^{\infty }{{{\left( -1 \right)}^{n}}\frac{{{\left( x-b \right)}^{n}}}{n!}{{p}_{n2}}-y\left( b \right)} \\ 
	& =-y\left( b \right)-\left( x-b \right){{p}_{12}}+\frac{{{\left( x-b \right)}^{2}}}{2!}{{p}_{22}}-\frac{{{\left( x-b \right)}^{3}}}{3!}{{p}_{32}}+\cdots  \\ 
	& =-y\left( b \right)-\left( x-b \right)\left( -\frac{y\left( a \right)}{\sinh\left( x-a \right)} \right)-\frac{{{\left( x-b \right)}^{3}}}{3!}\left( -\frac{y\left( a \right)}{\sinh\left( x-a \right)} \right)+\cdots  \\ 
	& =-y\left( b \right)+\frac{y\left( a \right)}{\sinh\left( x-a \right)}\left( \left( x-b \right)+\frac{{{\left( x-b \right)}^{3}}}{3!}+\frac{{{\left( x-b \right)}^{5}}}{5!}+\cdots  \right) \\ 
	& =-y\left( b \right)+\frac{y\left( a \right)}{\sinh\left( x-a \right)}\sinh\left( x-b \right)  
\end{aligned}	
\end{equation}

Therefore,
\begin{equation}\label{eq4.23}
\begin{aligned}
	y&=-\frac{{{Q}_{02}}}{{{Q}_{01}}}=-\frac{-y\left( b \right)+\frac{y\left( a \right)}{\sinh \left( x-a \right)}\sinh \left( x-b \right)}{\cosh \left( x-b \right)-\coth \left( x-a \right)\sinh \left( x-b \right)} \\ 
	& =-\frac{-y\left( b \right)\left( {{e}^{x-a}}-{{e}^{-\left( x-a \right)}} \right)+y\left( a \right)\left( {{e}^{x-b}}-{{e}^{-\left( x-b \right)}} \right)}{{{e}^{b-a}}-{{e}^{-\left( b-a \right)}}} \\ 
	& =\frac{y\left( a \right){{e}^{-b}}-y\left( b \right){{e}^{-a}}}{{{e}^{a-b}}-{{e}^{b-a}}}{{e}^{x}}+\frac{y\left( b \right){{e}^{a}}-y\left( a \right){{e}^{b}}}{{{e}^{a-b}}-{{e}^{b-a}}}{{e}^{-x}} 
\end{aligned}
\end{equation}

If using traditional method, the general solution is as follows:
\begin{equation}\label{eq4.24}
y\left( x \right)={{C}_{1}}{{e}^{x}}+{{C}_{2}}{{e}^{-x}}
\end{equation}

Substituting the function value conditions $x=a$, $y=y\left( a \right)$ and $x=b$, $y=y\left( b \right)$ into \eqref{eq4.24}, we get
\begin{equation}\label{eq4.25}
y\left( x \right)=\frac{y\left( a \right){{e}^{-b}}-y\left( b \right){{e}^{-a}}}{{{e}^{a-b}}-{{e}^{b-a}}}{{e}^{x}}+\frac{y\left( b \right){{e}^{a}}-y\left( a \right){{e}^{b}}}{{{e}^{a-b}}-{{e}^{b-a}}}{{e}^{-x}}
\end{equation} 	

The above formula is the same as \eqref{eq4.23}, i.e., the same solution is obtained by the traditional method, which shows the correctness of the method proposed in this paper. 

Although this method is more complex than the traditional method for simple problems, this method is a general method, and the procedure of solving process has consistency, not only for linear ordinary differential equations with constant coefficients, but also for linear ordinary differential equations with variable coefficients (see Example\ref{exm4.2}). 

This example illustrates the method proposed in this paper in detail. It can be seen that this method only needs algebraic computations (addition, subtraction, multiplication and division) and derivative calculations, does not need to solve the characteristic equations, does not need to determine the coefficients in the power series solution, and does not need integral calculations. 

\subsection{An example for sovling second order linear ordinary differential equation with variable coefficients}\label{subsec4.2}

\begin{Example}{Solve the following equation.}\label{exm4.2}
\begin{equation}\label{eq4.26}
	\left\{ \begin{aligned}
		& {y}''-xy=0 \\ 
		& {{\left. y \right|}_{x=b}}=y\left( b \right) \\ 
		& {{\left. y \right|}_{x=a}}=y\left( a \right) \\ 
	\end{aligned} \right.
\end{equation}
\end{Example}

\textbf{Solution}. Since
\begin{equation}\label{eq4.27}
\left\{ \begin{aligned}
	& {{p}_{21}}=0 \\ 
	& {{p}_{22}}=x \\ 
	& {{p}_{23}}=0 \\ 
\end{aligned} \right.
\end{equation}

From equation \eqref{eq4.1a}, we have
\begin{equation}\label{eq4.28}
{{p}_{31}}={{{p}'}_{21}}+{{p}_{22}}={{p}_{22}}=x
\end{equation}
\begin{equation}\label{eq4.29}
{{p}_{32}}={{{p}'}_{22}}=1
\end{equation}
\begin{equation}\label{eq4.30}
{{p}_{33}}=0
\end{equation}
\begin{equation}\label{eq4.31}
{{p}_{41}}=2{{{p}'}_{22}}=2
\end{equation}
\begin{equation}\label{eq4.32}
{{p}_{42}}={{{p}''}_{22}}+p_{22}^{2}={{x}^{2}}
\end{equation}
\begin{equation}\label{eq4.33}
{{p}_{43}}=0
\end{equation}
\begin{equation}\label{eq4.34}
{{p}_{51}}=3{{{p}''}_{22}}+p_{22}^{2}={{x}^{2}}
\end{equation}
\begin{equation}\label{eq4.35}
{{p}_{52}}={{{p}'''}_{22}}+4{{{p}'}_{22}}{{p}_{22}}=4x
\end{equation}
\begin{equation}\label{eq4.36}
{{p}_{53}}=0
\end{equation}
\begin{equation}\label{eq4.37}
{{p}_{61}}=4{{{p}'''}_{22}}+6{{{p}'}_{22}}{{p}_{22}}=6x
\end{equation}
\begin{equation}\label{eq4.38}
{{p}_{62}}=p_{22}^{\left( 4 \right)}+7{{{p}''}_{22}}{{p}_{22}}+4{{\left( {{{{p}'}}_{22}} \right)}^{2}}+p_{22}^{3}={{x}^{3}}+4
\end{equation}
\begin{equation}\label{eq4.39}
{{p}_{63}}=0
\end{equation}

From formulas \eqref{eq3.20}--\eqref{eq3.22}, let $k=2$, ${{x}_{1}}=b$, ${{x}_{2}}=a$,we can write the first six items of ${{Q}_{11}}$, ${{Q}_{12}}$, ${{Q}_{13}}$ in the following: 

\begin{equation}\label{eq4.40}
\begin{aligned}
	{{Q}_{11}}& =\sum\limits_{n=2}^{\infty }{{{\left( -1 \right)}^{n}}\frac{{{\left( x-a \right)}^{n}}}{n!}{{p}_{n1}}}-\left( x-a \right) \\ 
	& =-\left( x-a \right)-\frac{x{{\left( x-a \right)}^{3}}}{3!}+\frac{2{{\left( x-a \right)}^{4}}}{4!} -\frac{{{x}^{2}}{{\left( x-a \right)}^{5}}}{5!}+\frac{6x{{\left( x-a \right)}^{6}}}{6!}+\cdots   
\end{aligned}
\end{equation}
\begin{equation}\label{eq4.41}
\begin{aligned}
	{{Q}_{12}}& =\sum\limits_{n=2}^{\infty }{{{\left( -1 \right)}^{n}}\frac{{{\left( x-a \right)}^{n}}}{n!}{{p}_{n2}}}+1 \\ 
	& =1+\frac{x{{\left( x-a \right)}^{2}}}{2!}-\frac{{{\left( x-a \right)}^{3}}}{3!}+\frac{{{x}^{2}}{{\left( x-a \right)}^{4}}}{4!} -\frac{4x{{\left( x-a \right)}^{5}}}{5!}+\frac{\left( {{x}^{3}}+4 \right){{\left( x-a \right)}^{6}}}{6!}+\cdots   
\end{aligned}
\end{equation}
\begin{equation}\label{eq4.42}
{{Q}_{13}}=\sum\limits_{n=2}^{\infty }{{{\left( -1 \right)}^{n}}\frac{{{\left( x-a \right)}^{n}}}{n!}}{{p}_{n3}}-y\left( a \right)=-y\left( a \right)
\end{equation}

After descending order, we have
\begin{equation}\label{eq4.43}
\left\{ \begin{aligned}
	& {y}'={{p}_{11}}y+{{p}_{12}} \\ 
	& {{\left. y \right|}_{x=b}}=y\left( b \right) \\ 
\end{aligned} \right.
\end{equation}
where,
\begin{equation}\label{eq4.44}
\begin{aligned}
	{{p}_{11}}= -\frac{{{Q}_{12}}}{{{Q}_{11}}} 
	=-\frac{1+\frac{x{{\left( x-a \right)}^{2}}}{2!}-\frac{{{\left( x-a \right)}^{3}}}{3!}+\frac{{{x}^{2}}{{\left( x-a \right)}^{4}}}{4!}-\frac{4x{{\left( x-a \right)}^{5}}}{5!}+\frac{\left( {{x}^{3}}+4 \right){{\left( x-a \right)}^{6}}}{6!}+\cdots }{-\left( x-a \right)-\frac{x{{\left( x-a \right)}^{3}}}{3!}+\frac{2{{\left( x-a \right)}^{4}}}{4!}-\frac{{{x}^{2}}{{\left( x-a \right)}^{5}}}{5!}+\frac{6x{{\left( x-a \right)}^{6}}}{6!}+\cdots }  
\end{aligned}
\end{equation}
\begin{equation}\label{eq4.45}
\begin{aligned}
	{{p}_{12}}=-\frac{{{Q}_{13}}}{{{Q}_{11}}} =\frac{y\left( a \right)}{-\left( x-a \right)-\frac{x{{\left( x-a \right)}^{3}}}{3!}+\frac{2{{\left( x-a \right)}^{4}}}{4!}-\frac{{{x}^{2}}{{\left( x-a \right)}^{5}}}{5!}+\frac{6x{{\left( x-a \right)}^{6}}}{6!}+\cdots }  
\end{aligned}
\end{equation}

From \eqref{eq3.33}--\eqref{eq3.34}, let $k=2$, ${{x}_{1}}=b$, ${{x}_{2}}=a$, we have
\begin{equation}\label{eq4.46}
\begin{aligned}
	{{Q}_{01}}& =\sum\limits_{n=1}^{\infty }{{{\left( -1 \right)}^{n}}\frac{{{\left( x-b \right)}^{n}}}{n!}{{p}_{n1}}+1} \\ 
	& =1-\left( x-b \right){{p}_{11}}+\frac{{{\left( x-b \right)}^{2}}}{2!}{{p}_{21}}-\frac{{{\left( x-b \right)}^{3}}}{3!}{{p}_{31}}+\frac{{{\left( x-b \right)}^{4}}}{4!}{{p}_{41}}-\cdots   
\end{aligned}
\end{equation}
\begin{equation}\label{eq4.47}
\begin{aligned}
	{{Q}_{02}}& =\sum\limits_{n=1}^{\infty }{{{\left( -1 \right)}^{n}}\frac{{{\left( x-b \right)}^{n}}}{n!}{{p}_{n2}}-y\left( b \right)} \\ 
	& =-y\left( b \right)-\left( x-b \right){{p}_{12}}+\frac{{{\left( x-b \right)}^{2}}}{2!}{{p}_{22}}-\frac{{{\left( x-b \right)}^{3}}}{3!}{{p}_{32}}+\cdots   
\end{aligned}
\end{equation}

So the equation \eqref{eq4.26} is solved in the following:
\begin{equation}\label{eq4.48}
y=-\frac{{{Q}_{02}}}{{{Q}_{01}}}=\frac{\sum\limits_{n=1}^{\infty }{{{\left( -1 \right)}^{n}}\frac{{{\left( x-b \right)}^{n}}}{n!}{{p}_{n2}}-y\left( b \right)}}{\sum\limits_{n=1}^{\infty }{{{\left( -1 \right)}^{n}}\frac{{{\left( x-b \right)}^{n}}}{n!}{{p}_{n1}}+1}}
\end{equation}
where ${{p}_{n1}}$ and ${{p}_{n2}}$ are obtained by the following recursive formulas $($ let $k=1$ in \eqref{eq3.13} $)$:
\begin{equation}\label{eq4.49}
\left\{ \begin{array}{ll}
	{{p}_{n1}}={{{{p}'}}_{\left( n-1 \right)1}}+{{p}_{\left( n-1 \right)1}}{{p}_{11}} \\ 
	{{p}_{n2}}={{{{p}'}}_{\left( n-1 \right)2}}+{{p}_{\left( n-1 \right)1}}{{p}_{12}} \\ 
\end{array}
\right.
,
\begin{array}{l}
	n=2,3,\cdots
\end{array}
\end{equation}

\section{Conclusions and prospects}\label{sec5}
This section summarizes the characteristics and advantages of the results and gives prospects for the future.
\begin{itemize}
\item  All the theories and methods in this paper are my original innovations. All the theories and methods proposed in this paper are different from the traditional known methods.
\end{itemize}

\begin{itemize}
\item It is generally considered that the solution of linear ordinary differential equations of order $k$ ($k\ge 2$) is very difficult to obtain. In this paper, a general method for solving linear ordinary differential equations of order $k$ ($k\ge 2$) is given.
\end{itemize}

\begin{itemize}
\item The method proposed in this paper is a general method, the steps of solving the equation are unified and consistent, it can be used to solve the linear ordinary differential equations, not only for solving the linear ordinary differential equations with constant coefficients, but also for solving the linear ordinary differential equations with variable coefficients.
\end{itemize}

\begin{itemize}
\item The core process is recursion and reduction of order.
\end{itemize}

\begin{itemize}
\item Only algebraic computation (addition, subtraction, multiplication and division) and derivative calculation are required.
\end{itemize}

\begin{itemize}
\item The function value conditions are automatically satisfied.
\end{itemize}

\begin{itemize}
\item For the function value conditions, a general method for solving linear ordinary differential equations is given in this paper. For the conditions of derivative value, I will study it in one of my follow up papers.
\end{itemize}

The following are some questions that need further study, some of research results have been made, but some problems still need to be further explored.

\begin{itemize}
\item The solution of the function value problem of order $k$ \eqref{eq3.3} is given in this paper. For a linear ordinary differential equation with other conditions, how do we give the result?
\end{itemize}

\begin{itemize}
\item Which type of functions does the coefficient ${{p}_{ki}}$ ($i=1,2,\cdots,k+1$) in equation \eqref{eq3.3} satisfy condition \eqref{eq2.13}? How can we get some sufficient conditions? 
\end{itemize}

\begin{itemize}
\item The solution of the function value problem of order $k$ of differential equation \eqref{eq3.3} is a special solution. If we use the formula \eqref{eq2.8} containing an arbitrary constant to construct differential equations according to the method proposed in this paper, the solution obtained will contain $k$ arbitrary constants. The question is, is this solution containing $k$ arbitrary constants the general solution of differential equation \eqref{eq3.1}? If so, how to prove it?
\end{itemize}

\begin{itemize}
\item How to apply the results proposed in this paper to numerical analysis and obtain approximate solutions of differential equations?
\end{itemize}

\begin{itemize}
\item Set ${{Z}^{\infty }}\left( \mathbf{I} \right)$ is the set of all the solutions of differential equation \eqref{eq2.12}, so it is very important for the study of set ${{Z}^{\infty }}\left( \mathbf{I} \right)$. Is the set ${{Z}^{\infty }}\left( \mathbf{I} \right)$ a linear space? If so, what is its base and dimension?
\end{itemize}

\begin{itemize}
\item What is the internal relationship between ${{Z}^{\infty }}\left( \mathbf{I} \right)$ and ${{C}^{\infty }}\left( \mathbf{I} \right)$?
\end{itemize}

There are still many questions that are hard to put forward in one paper. For the new theory and method put forward in this paper, a series of problems can be formed, and it is impossible to solve all the problems in one paper. This is my first paper on this topic, so please allow me to name the title of the paper as "General methods for solving ordinary differential equations 1". I hope there will be a "General methods for solving ordinary differential equation 2", a "General methods of solving ordinary differential equation 3", and so on. I will continue to study this topic seriously. If my research results can arouse the interest of fellow scholars, it will be a great honor for me.


\end{document}